\documentclass{article}
\usepackage{amsfonts,amssymb,amsmath}
\usepackage{graphicx}
\newcommand{\R}{\mathbb{R}}

\begin{document}
\title{Ultrametric Embedding: Application to Data Fingerprinting and to Fast
Data Clustering}

\author{Fionn Murtagh \\
Department of Computer Science, \\
Royal Holloway, University of London, \\
Egham, Surrey TW20 0EX, England. \\ 
fmurtagh@acm.org}

\maketitle

\begin{abstract}
We begin with pervasive ultrametricity due to high dimensionality and/or
spatial sparsity.  How extent or degree of ultrametricity can be quantified
leads us to the discussion of varied practical cases when ultrametricity can 
be partially or locally present in data.  We show how the ultrametricity 
can be assessed in text or document collections, and 
in time series signals.   An aspect of importance here is that to draw benefit 
from this perspective the data may need to be recoded.  Such data recoding can also 
be powerful in proximity searching, as we will show, where the data is embedded globally 
and not locally in an ultrametric space.
\end{abstract}

%\begin{keyword}      
%multivariate data analysis, cluster analysis, hierarchy, 
%factor 
%analysis, correspondence analysis, ultrametric, p-adic, phylogeny
%\end{keyword}

\section{Introduction}

The topology or inherent shape and form of an object is important.  In 
data analysis, the inherent form and structure of data clouds are important.
Quite a few models of data form and structure are used in data analysis.  
One of them is a hierarchically embedded set of clusters, -- a 
hierarchy.  It is traditional (since at least the 1960s) to impose such a
form on data, and if useful to assess the goodness of fit.  Rather than
fitting a hierarchical structure to data, our recent work
has taken a different orientation: we seek to find (partial or global) 
inherent hierarchical structure in data.  As we will describe in this 
article, there are interesting findings that result from this, and some
very interesting perspectives are opened up for data analysis.  

A formal definition of hierarchical structure is provided by ultrametric
topology (in turn, related closely  to p-adic number theory).  We will
return to this in section \ref{sect23} below.  First, though, we will 
summarize some of our findings.

Ultrametricity is a pervasive property of observational data.  It 
arises as a limit case when data dimensionality or sparsity grows.  More
strictly such a limit case is a regular lattice structure and
ultrametricity is one possible representation for it.  Notwithstanding 
alternative representations, ultrametricity offers
 computational efficiency (related to tree depth/height being logarithmic 
in number of terminal nodes), linkage with dynamical or related 
functional properties (phylogenetic interpretation), and 
processing tools based on well known p-adic or ultrametric theory (examples: 
deriving a partition, or applying an ultrametric wavelet transform).  

Local ultrametricity is also of importance.  
Practical data sets (derived from, or observed in,
 databases and data spaces) present some but not exclusively 
ultrametric characteristics.  This can be used for forensic data 
exploration (fingerprinting data sets, as we discuss below).  
Local ultrametricity has been used to expedite search and discovery in 
information spaces (in \cite{cha2} as discussed by us in \cite{ref08}, which we 
will not discuss further here).  Such proximity searching and matching has traditionally
been addressed ultrametrically by fitting a hierarchy to data.  Below, we show a 
different way to embed the data (in a computationally highly efficient way)
in an ultrametric space, using a principle employed
in our local ultrametric work: namely, data recoding.  

Our ultimate aim in this work is  to proceed 
a lot further, and gain new insights into data (and 
observed phenomena and events) through ultrametric (topology) or equivalently p-adic 
(algebra) representation theory.

\section{Quantifying Degree of Ultrametricity}
\label{sect23}

Summarizing a full description in Murtagh \cite{ref08} we explored two
measures quantifying how ultrametric a data set is, -- Lerman's and a new
approach based on triangle invariance (respectively, 
the second and third approaches described in this section).

The triangular inequality holds for a metric space: $d(x,z) \leq 
d(x,y) + d(y,z)$ for any triplet 
of points $x,y,z$.  In addition the properties 
of symmetry and positive definiteness are respected.  The ``strong 
triangular inequality'' or ultrametric inequality is: $d(x,z) \leq 
\mbox{ max } \{ d(x,y), d(y,z) \}$ for any triplet $x,y,z$.  An
ultrametric
space implies respect for a range of stringent properties.  For
example, 
the triangle formed by any triplet is necessarily isosceles, with the
two
large sides equal; or is equilateral.

\begin{itemize}
\item 
Firstly, Rammal et al.\ \cite{ref09b} used discrepancy between each pairwise 
distance and the corresponding subdominant ultrametric.  Now, the
subdominant ultrametric is also known as the ultrametric distance
resulting from the single linkage agglomerative hierarchical
clustering method.   Closely related graph structures include the
minimal spanning tree, and graph (connected) components.  
While the subdominant provides a good fit to the
given distance (or indeed dissimilarity), it suffers from the
``friends of friends'' or chaining effect.  

\item
Secondly, Lerman \cite{ref07} developed a measure of ultrametricity,
termed H-classifiability,  using
ranks of all pairwise given distances (or dissimilarities).  The
isosceles (with small base) or equilateral requirements of the
ultrametric inequality impose constraints on the ranks.  The interval between
median and maximum rank of every set of triplets must be empty for 
ultrametricity.  We have used extensively  
Lerman's measure of degree of ultrametricity in a data set.  
Taking ranks provides scale invariance.  
But the limitation of Lerman's approach, 
we find, is that it is not reasonable to
study ranks of real-valued (values in non-negative reals) 
distances defined on a large set of points.

\item
Thirdly, our own measure of extent of ultrametricity \cite{ref08}
can be described algorithmically.  We 
examine triplets of points (exhaustively if possible, or otherwise
through sampling), and determine the three angles formed by the
associated triangle.  We select the smallest angle formed by the triplet
points.  Then we check if the other two remaining angles are
approximately equal.  If they are equal then our triangle is isosceles
with small base, or equilateral (when all triangles are equal).  The
approximation to equality is given by 2 degrees (0.0349 radians).  
Our motivation for
the approximate (``fuzzy'') equality is that it makes our approach
robust and independent of measurement precision.  
\end{itemize}

A supposition for use of our measure of ultrametricity is that we can can 
define angles (and hence triangle properties).  This in turn presupposes a
scalar product.  Thus we presuppose a normed vector space with a scalar product --
a Hilbert space -- to provide our needed environment.  Quite a general way to 
embed data, to be analyzed, in a Euclidean space, is to use correspondence 
analysis \cite{ref08888}.  This explains our interest in using correspondence 
analysis quite often in this work: it provides a convenient and versatile way to 
take input data in many varied formats (e.g., ranks or scores, presence/absence,
frequency of occurrence, and many other forms of data) and map them into a 
Euclidean, factor space.

\section{Ultrametricity and Dimensionality}

\subsection{Distance Properties in Very Sparse Spaces}
\label{sect22}

Murtagh \cite{ref08}, and earlier work by Rammal et al.\ \cite{ref09a,ref09b},
has demonstrated the pervasiveness of ultrametricity, by 
focusing on the fact that  
sparse high-dimensional data tend to be ultrametric.  
In Murtagh \cite{ref08} it is shown how numbers of points
in our clouds of data points are irrelevant; but what counts is the
ambient spatial dimensionality.  Among cases looked at are statistically 
uniformly
(hence ``unclustered'', or without structure in a certain sense)
distributed points, and statistically 
uniformly distributed hypercube vertices (so the
latter are random 0/1 valued vectors).  Using our ultrametricity
measure, there is a clear tendency to ultrametricity as the spatial
dimensionality (hence spatial sparseness) increases.  

%Figure \ref{fig6} illustrates these findings.  Dimensionality is increased
%up to 5000.   The number of points is kept at 5000, but varying this causes 
%little if any difference in the results.  (Exponentially increasing this 
%number of points, $n$, is quite a different matter, but is of no interest
%to us here.)  The results are quite similar irrespective of the 
%random number generation seed used, and also with respect to different 
%random generations of the data.  
%Two different sets of experiments are run, and the results 
%are displayed in the two curves.  Firstly, uniformly distributed real 
%values are
%used (hence $n$ points, in $\R^m$, where $m$ is the dimensionality), and 
%secondly uniformly distributed hypercube vertices are generated (hence
%$n$ points, in $ \{ 0, 1 \}^m$).  Using our triangle-based quantification of 
%ultrametricity, we find high dimensionality to approach the limit of 
%complete ultrametricity.  

%\begin{figure}
%\centering
%\includegraphics[width=9cm,angle=270]{fig-6-ultrametric.ps}
%\caption{Upper curve: uniformly distributed values.  Lower curve:
%random hypercube vertex points. A value of our triangle-based 
%ultrametricity measure equal to 1 is related
%to global ultrametricity.  For each dimensionality (50, 100, 500,
%1000, 5000) we used number of points $n = 5000$.  In other experiments
%we found very little variation as a function of $n$.}
%\label{fig6}
%\end{figure}

As \cite{hall} also show, Gaussian data behave in the same way and a demonstration
of this is seen in Table \ref{tabunifgauss}.  To provide an idea of consensus of these
results, the 200,000-dimensional Gaussian was repeated and yielded on successive runs
values of the ultrametricity measure of: 0.96, 0.98, 0.96.  

\begin{table}
\begin{center}
\begin{tabular}{lllll} \hline
No. points &    Dimen. &  Isosc. &   Equil. &    UM   \\ \hline    
           &           &         &          &         \\
Uniform    &           &         &          &         \\
           &           &         &          &         \\
100        &    20     &    0.10 &    0.03  &    0.13 \\
100        &    200    &    0.16 &    0.20  &    0.36 \\
100        &    2000   &    0.01 &    0.83  &    0.84 \\
100        &    20000  &    0    &    0.94  &    0.94 \\
100        &    200000 &    0    &    0.97  &    0.97 \\
           &           &         &          &         \\
Hypercube  &           &         &          &         \\
           &           &         &          &         \\
100        &    20     &    0.14 &   0.02   &    0.16 \\
100        &    200    &    0.16 &   0.21   &    0.36 \\
100        &    2000   &    0.01 &   0.86   &    0.87 \\
100        &    20000  &    0    &   0.96   &    0.96 \\
100        &    200000 &    0    &   0.97   &    0.97 \\  
           &           &         &          &         \\
Gaussian   &           &         &          &         \\
           &           &         &          &         \\
100        &    20     &    0.12 &    0.01  &    0.13 \\
100        &    200    &    0.23 &    0.14  &    0.36 \\
100        &    2000   &    0.04 &    0.77  &    0.80 \\
100        &    20000  &    0    &    0.98  &    0.98 \\
100        &    200000 &    0    &    0.96  &    0.96 \\ \hline
\end{tabular}
\end{center}
\caption{Typical results, based on 300 sampled triangles from triplets of
points.  For uniform, the data are generated on [0, 1]; hypercube vertices
are in $\{ 0, 1\}^{\mbox{Dimen}}$, 
and for Gaussian, the data are
of mean 0, and variance 1.  Dimen. is the ambient
dimensionality.  Isosc. is the number of isosceles triangles with small base, as a
proportion of all triangles sampled.  Equil. is the number of equilateral triangles
as a proportion of triangles sampled.  UM is the proportion of ultrametricity-respecting
triangles (= 1 for all ultrametric).}
\label{tabunifgauss}
\end{table}

In the following, we explain why high dimensional and/or sparsely populated
 spaces are ultrametric.  

As dimensionality grows, so too do distances (or indeed
dissimilarities, if
they do not satisfy the triangular inequality).  The least change
possible for dissimilarities to become distances has been formulated
in terms of the smallest additive constant needed, to be added to all
dissimilarities \cite{tor,cai1,cai2,neu}.
Adding a sufficiently large 
constant to all dissimilarities transforms them into a set of
distances.  Through addition of a larger constant, it follows that
distances become approximately equal, thus verifying a trivial case of
the ultrametric or ``strong triangular'' inequality.  Adding to
dissimilarities or distances may be a direct consequence of increased
dimensionality. 

For a close fit or good approximation,  
the situation is not as simple for taking dissimilarities, or
distances,
into ultrametric distances.  A best fit solution is given by \cite{desoete}
(and software is available in R \cite{hornik}).
If we want a close fit to the given 
dissimilarities then a good choice would avail either of the maximal 
inferior, or subdominant, ultrametric; or the minimal superior
ultrametric.
Stepwise algorithms for these are commonly known as, respectively,
single
linkage hierarchical clustering; and complete link hierarchical
clustering.  
(See \cite{ref2,ref07,mur85a} and other texts on 
hierarchical clustering.)  

\subsection{No ``Curse of Dimensionality'' in Very High Dimensions}
\label{sect21}

Bellman's \cite{bel}   ``curse of dimensionality'' relates to exponential 
growth of hypervolume as a function of dimensionality.  Problems
become tougher as dimensionality increases.  In particular problems 
related to proximity search in high-dimensional spaces tend to become
intractable.  

In a way, a ``trivial limit'' (Treves \cite{ref120}) 
case is reached as dimensionality increases.  This makes high
dimensional 
proximity search very different, and given an appropriate data
structure -- such as a binary hierarchical clustering tree -- we can
find nearest neighbors in worst case $O(1)$ or constant computational
time \cite{ref08}.  The proof is simple: the tree data structure
affords a constant number of edge traversals.  

The fact that limit properties are ``trivial'' makes them no less 
interesting to study.  Let us refer to such ``trivial'' properties as
(structural or geometrical) regularity properties (e.g.\ all points
lie on a regular lattice).  First of all, the symmetries of regular
structures in our data may be of importance.  Secondly, ``islands'' or
clusters in our data, where each ``island'' is of regular structure,
may be exploitable.  Thirdly, the mention of exploitability points to
the application areas targeted: in this article, we focus on search and
matching and show some ways in which ultrametric regularity can be
exploited in practice.  Fourthly, and finally, regularity by no means
implies complete coverage (e.g., existence of all pairwise linkages)
implying that interesting or revealing structure will be present in 
observed or recorded data sets. 

\begin{figure}
\includegraphics[width=16cm]{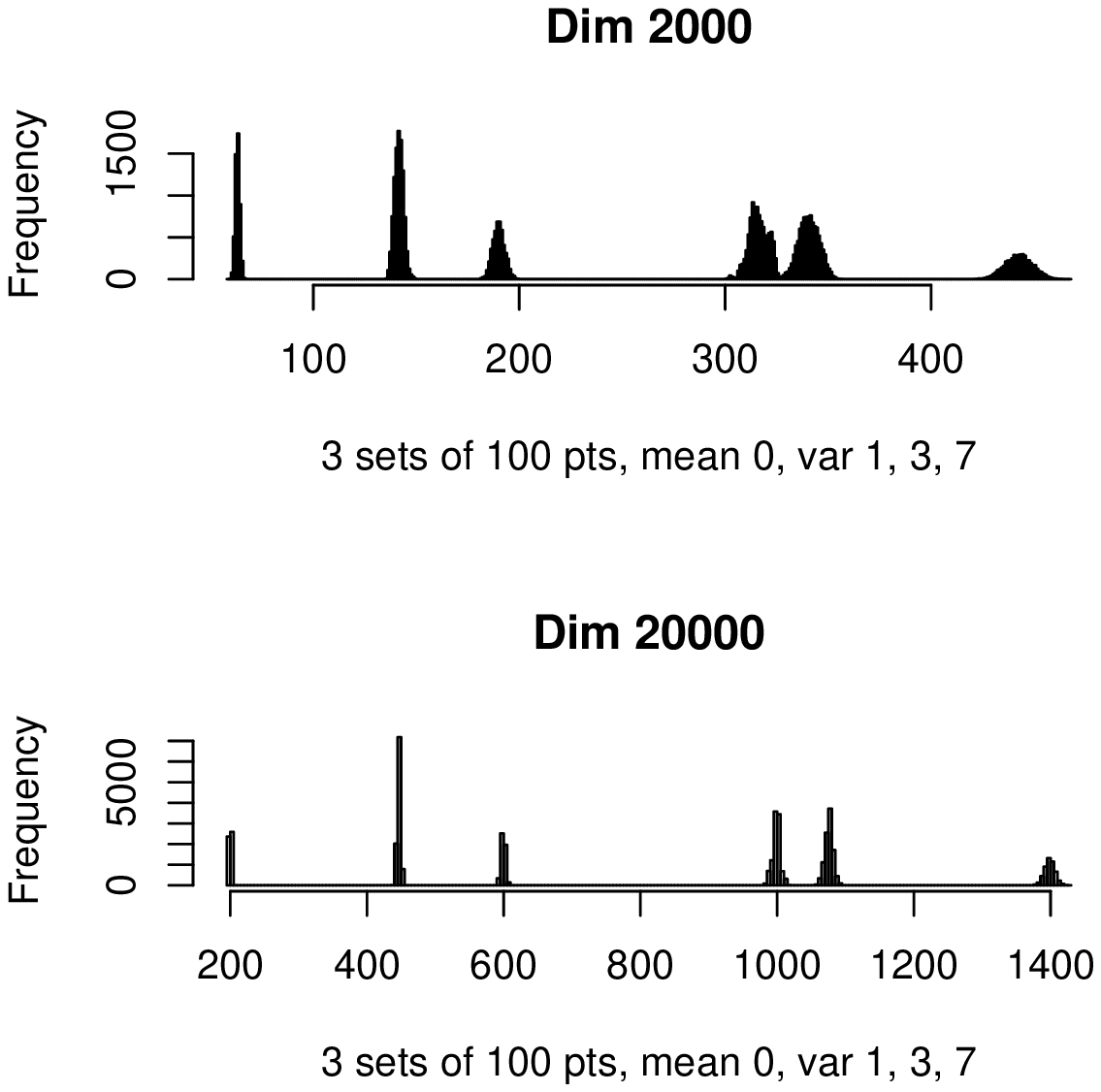}
\caption{An illustration of how ``symmetry'' or ``structure'' can become 
increasingly pronounced as dimensionality increases.
Shown are two simulations, each with 3 sub-populations of Gaussian-distributed data,
in, respectively, ambient dimensions of 2000 and 20,000.}
\label{fig1}
\end{figure}

Thus we see that in very high dimensions, and/or in very (spatially) sparse data clouds, 
there is a simplification of structure, 
which can be used to mitigate any ``curse of dimensionality''.  Figure \ref{fig1}
shows how the distances within and between clusters become tighter with increase in 
dimensionality.  

\section{Increasing Ultrametricity Through Data Recoding}

\subsection{Ultrametricity of Text}
\label{sect42}

In \cite{ref0888}, words appearing in a text 
(in principle all, but in practice a set of the few
hundred most frequent)  are used to fingerprint  the text.  
Rare words in a text corpus may be appropriate for querying the 
corpus for relevant texts, but such words are of little help for 
inter-text characterization and comparison.  We also use entire words, with 
no stemming or other preprocessing.  
A full justification for such an approach to textual 
data analysis can be found in Murtagh \cite{ref08888}.  

So our methodology for studying a set of texts is to characterize each 
text with numbers of terms appearing in the text, for a set of terms.  
The $\chi^2$ distance
is an appropriate weighted Euclidean distance for use with such data
\cite{ref2,ref0888}.  
Consider texts $i$ and $i'$ crossed by words $j$.  Let $k_{ij}$ be the 
number of
occurrences of word $j$ in text $i$.  Then, omitting a constant, 
the $\chi^2$ distance between texts $i$ and $i'$ is given by 
$ \sum_j 1/k_j ( k_{ij}/k_i - k_{i'j}/k_{i'} )^2$.  The weighting term is 
$1/k_j$.  The weighted Euclidean distance is between the {\em profile} 
of text $i$, viz.\ $k_{ij}/k_i$ for all $j$, and the analogous 
{\em profile} of text $i'$.  (Our discussion is to within a constant because
we actually work on {\em frequencies} defined from the numbers of 
occurrences.) 

Correspondence analysis allows us to project the space of documents (we could
equally well explore the terms in the {\em same} projected space) into 
a Euclidean space.  It maps the all-pairs $\chi^2$ distance into the 
corresponding Euclidean distance.  In the resulting factor space, we use
our triangle-based approach for quantifying how ultrametric the data are.

We did this, \cite{ref0888},  
for a large number of texts (3 Jane Austen novels, James 
Joyce's {\em Ulysses}, 
technical reports -- 50 airline accident reports from the NTSB, National 
Transport Safety Board, fairy tales -- 209 fables of the Brothers
Grimm, 214 dream reports from the DreamBank repository, 
Aristotle's {\em Categories}, etc.), 
finding consistent degree of 
ultrametricity results over texts of the same sort.

Some very intriguing ultrametricity characterizations were found in our
work.  For example, we found that the technical vocabulary of air accidents 
did not differ greatly in terms of inherent ultrametricity compared to the 
Brothers Grimm fairy tales.  Secondly we found that novelist Austen's 
works were distinguishable from the Grimm fairy tales.  Thirdly we found 
dream reports to be have higher ultrametricity level than the other 
text collections.   Dream-like ultrametric characteristics of Joyce's
{\em Ulysses} were also apparent.

Given that  local ultrametricity on sets of terms or words implies locally
hierarchical relationships between them, we have pursued this work 
in the direction of automating the task of ontology construction.  See 
\cite{autoonto}.  

\subsection{Data Recoding in the Correspondence Analysis Tradition}

If the $\chi^2$ distance (see above, section \ref{sect42}) 
is used on data tables with constant
marginal sums then it becomes a weighted Euclidean distance.  This is 
important for us, because it means that we can directly influence the 
analysis by equi-weighting, say, the table rows in the following way: 
we double the row vector values by including an absence (0 value) 
whenever there is a presence (1 value) and vice versa.  Or for a table of 
percentages, we take both the original value $x$ and $100 - x$.  In the 
correspondence analysis tradition \cite{ref2,ref088} this is known as 
{\em doubling} ({\em d\'edoublement}).  

More generally, booleanizing, or making qualitative, data in this way, for a
varying (value-dependent) number of target value categories (or modalities) 
leads to the form of coding known as {\em complete disjunctive form}.  

Such coding increases the embedding dimension, and data sparseness, and thus
may encourage degree of ultrametricity.  That it can do more we will now 
show.  

The iris data has been very widely used as a toy data set since 
Fisher used it in 1936 (\cite{fish}, taking from a 1935 article by Anderson) 
to exemplify discriminant analysis.  It consists of 150 iris flowers, 
each characterized by 4 petal and sepal, width and breadth, measurements.  
On the one hand, therefore, we have the 150 irises in $\R^4$.  Next, 
each variable value was recoded to be a rank (all ranks of a given variable
considered) and the rank was boolean-coded (viz., for the top rank variable 
value, $1000 \dots$, for the second rank variable value, $0100 \dots$, etc.).
Following removal of zero total columns, the second data set 
defined the 150 irises in $\R^{123}$.  Actually, this definition of the 
150 irises is in fact in $\{0,1\}^{123}$.  

Our triangle-based measure of the degree of ultrametricity in a data set
(here the set of irises), with 0 = no ultrametricity, and 1 = every triangle
an ultrametric-respecting one, gave the following: for irises in 
$\R^4$, 0.017; and for irises in $\{0,1\}^{123}$: 0.948.  
% See G4, UMexpts 

This provides a nice illustration of how recoding can dramatically change 
the picture provided by one's data.  Furthermore it provides justification 
for data recoding if the ultrametricity can be instrumentalized by us 
in some way, e.g.\ to facilitate fast proximity search.   

\subsection{Ultrametricity of Time Series}
\label{sect5353}

In Murtagh \cite{ref088} we use the following coding to show that chaotic time 
series are less ultrametric than, say, financial (futures, FTSE -- 
Financial Times Stock Exchange index, stock price index), biomedical (EEG for
normal and epileptic subjects, eyegaze trace), telecoms (web traffic)  or 
meteorological (Mississippi water level, sunspots) time 
series; random generated (uniformly distributed) time series data are 
remarkably similar in their ultrametric properties; and ultrametricity 
can be used to distinguish various types of biomedical (EEG) signals.  

A time series can be easily embedded in a space of dimensionality $m$, 
by taking successive intervals of length $m$, or a delay embedding of 
order $m$.  Thus we define points 
$$
\mathbf{x}_r = ( x_{r-m+1}, x_{r-m+2}, \dots , x_{r-1}, x_r )^t \in \R^m$$ 
where $t$ denotes vector transpose.  
%Based on previous results we expect that
%as the dimension $m$ grows, then the point set in $\R^m$ becomes more 
%ultrametric.  This finding is borne out below.

Given any $\mathbf{x}_r = ( x_{r-m+1}, x_{r-m+2}, \dots , x_{r-1}, x_r )^t \in 
\R^m$, let us consider the set of $s$ such contiguous intervals
determined from the time series of overall size $n$.  For 
convenience we will take $s = \lfloor n/m \rfloor$ where $\lfloor . \rfloor$
is integer truncation.  The contiguous intervals could be 
overlapping but for exhaustive or near-exhaustive coverage it is 
acceptable that they be non-overlapping.  In our work, the
intervals were non-overlapping.  
The quantification of the ultrametricity of the overall time series is 
provided by the aggregate over $s$ time intervals of the ultrametricity of 
each  $\mathbf{x}_r$, 
$1 \leq r \leq s$.  

We seek to directly quantify the extent of 
ultrametricity in time series data.  Earlier in this article we have seen 
how increase in ambient spatial 
dimensionality leads to greater ultrametricity.  However 
it is not satisfactory from a 
practical point of view to simply increase the embedding dimensionality $m$
insofar as short memory relationships are of greater practical 
relevance (especially for prediction).
The greatest possible 
value of $m > 1 $ is the total length of the time series, $n$.   
Instead we will look for an ultrametricity measurement approach for given 
and limited sized dimensionalities $m$.  Our experimental results
for real and for random data sets are for ``window'' lengths 
$m = 5, 10, \dots , 105, 110$.  

We seek local ultrametricity, i.e.\ hierarchical structure, by 
studying the following:
Euclidean distance squared, $d_{jj'} = (x_{rj} - x_{rj'})^2$ 
for all $ 1 \leq j, j' \leq m$ in each time 
window, $\mathbf{x}_r$.  It will be noted below in this section 
how this assumption of Euclidean distance squared has worked well but is 
not in itself important: in principle any dissimilarity can be used.  

We enforce sparseness \cite{ref09a,ref09b,ref08} on our given distance
values, $\{ d_{jj'} \}$.  We do this by thresholding
each unique value $d_{jj'}$, in the range $ \mbox{max}_{jj'} d_{jj'} - 
\mbox{min}_{jj'} d_{jj'} $, by an integer in $\{ 1, 2 \}$.  Note that the 
range is chosen with reference to the currently considered time series
window, $ 1 \leq j, j' \leq m$.  
%Note too that the value of 
%$p$ must be specified.  In our work we set  $p = 2$.  
Thus far, the recoded value, $d'_{jj'}$
is not necessarily a distance.  With the extra requirement that 
$d'_{jj'} \longrightarrow 0$ whenever $j = j'$ it can be shown that 
$d'_{jj'}$ is a metric \cite{ref088}.

To summarize, in our coding,
a small pairwise transition is mapped onto a value of
1; and a large pairwise transition is mapped onto a value of
2.  A pairwise transition is defined not just for data values that are
successive in time but for any pair of data values in the window considered.

This coding can be considered as (i) taking a local region, defined by 
the sliding window, and (ii) coding pairwise ``change'' = 2, versus 
``no change'' = 1, relationships.  Then, based on these new distances,
we use the ultrametric triangle properties to assess conformity to 
ultrametricity.   The average overall ultrametricity in the time series, 
quantified in this way, allows us to fingerprint our time series.

A wide range of window sizes (i.e., lengths), $m$, was investigated.
Window size is not important: in relative terms the results 
found remain the same.  Taking part of a time series and comparing the 
results to the full time series gave similar outcomes, thus indicating 
that the fingerprinting was an integral property of the data.  

Our ``change/no change'' metric is crucial here, and not the input 
dissimilarity which is mapped onto it.  Note too that generalization to 
multivariate time series is straightforward.  

Eyegaze
trace signals were found to be remarkably high in ultrametricity, which 
may be due to extreme values (truncated 
off-scale readings resulting from the subject's
blinking) that were not subject to preprocessing.  Web traffic was
also very high in ultrametricity, due to to extreme values.
All EEG data sets
were close together, with clear separation between the normal sleep 
subject, and the epilepsy cases.  The lowest ultrametricity 
was found for chaotic time series.    

\subsection{Fast Clustering through Baire Space Embedding}
\label{sectsearch}

The clustering of
chemical compounds, based on chemical descriptors or representations, 
is important in the pharmaceutical and chemical sectors.  It is used 
for screening and knowledge discovery in large databases of chemical 
compounds.  A chemical compound  is encoded (through various schemes that 
are not of relevance to us here) as a fixed length  bit string 
(i.e.\ a set of boolean or 0/1 values).  We have started to look at a 
set of 1.2 million chemical compounds, each characterized (in a given 
descriptor or coding system, the Digital Chemistry bci1052 dictionary of fragments) 
by 1052 variables.  

While attributes per chemical compound are roughly Gaussian in distribution,
chemicals per attribute follow a power law.  
We found the probability of having more than 
$p$ chemicals per attribute to be approximately $c/p^{1.23}$ for large $p$
and for constant, $c$.  This warrants normalization, which we do by dividing 
attribute/chemical presence values by the attribute marginal (i.e., attribute 
column sum).  Any presence value is now a floating point value.

Consider now the very simplified  example of two chemicals, $x$ and $y$,
with just one attribute, whose maximum precision of measurement is $K$. 
So let us consider  $x_K = 0.478$; and $y_K = 0.472$. 
In these cases, maximum precision, $| K | = 3$.  For first decimal place 
$k = 1$, we find $x_k
= y_k = 4$.  For $k = 2, x_k = y_k$.  But for $k = 3, x_k \neq y_k$.
We now introduce the following distance:
%\begin{eqnarray*}
%d_B(x_K, y_K) =    & 1                   & \mbox{if } x_1 \neq y_1  \\
%                   & \mbox{inf } 2^{-n}  & x_n = y_n  \ \ \  1 \leq n \leq |K|
%\end{eqnarray*}
\begin{equation*}
d_B(x_K, y_K) =
\begin{cases}
1                   & \text{if } x_1 \neq y_1  \\
\text{inf } 2^{-n}  & x_n = y_n  \ \ \  1 \leq n \leq |K|
\end{cases}
%\label{eqn2}
\end{equation*}
So here $d_B(x_K, y_K) = 2^{-3}$.  This distance is a greatest common prefix metric, and
indeed ultrametric.  Its maximum value is 1, i.e.\ it is a 1-bounded ultrametric.  Our 
reason for use of $d_B$ to denote this distance is due to it endowing a metric on 
the Baire space, the space of  countably infinite sequences.  

The case of multiple attributes is handled as follows.  We have the set $J$ of 
attributes.  Hence 
we have $|J|$ values for each chemical structure.  So the $i$th
chemical structure, for each $j \in J$ value with precision $|K|$, is $x_{iJK}$\
.
Collectively, all our data are expressed by $x_{IJK}$.  As before, we
normalize by column
sums to work therefore on $x_{IJK}^J$.  To find the Baire distance properties
we work simultaneously on all $J$ values, corresponding to a given
chemical structure.
Therefore the partition at level $k = 1$ has
clusters defined as all those numbers indexed by $i$ that share the same
$k = 1$, or 1st, digit {\em in all $J$ values}.

\begin{table}
\begin{center}
\begin{tabular}{|c|r|}\hline
Sig. dig. $k$ &  No. clusters  \\ \hline
           &                \\
  4        &    6591        \\
  4        &    6507        \\
  4        &    5735        \\
           &                \\
  3        &    6481        \\
  3        &    6402        \\
  3        &    5360        \\
           &                \\
  2        &    2519        \\
  2        &    2576        \\
  2        &    2135        \\
           &                \\
  1        &     138        \\
  1        &     148        \\
  1        &     167        \\
           &                \\ \hline
\end{tabular}
\end{center}
\caption{Results for the three different data sets, each consisting of
7500 chemicals,  are shown in immediate
succession.  The number of significant decimal digits is 4 (more precise,
and hence more different clusters found), 3, 2, and 1 (lowest precision
in terms of significant digits).}
\label{tab2}
\end{table}

Table \ref{tab2} demonstrates how this works.  
In Table \ref{tab4} we look at k-means, using as input the cluster centers
provided by the 1-significant digit Baire approach.  Relatively very few
changes were found.  We note that the partitions in each case are
dominated by a very large cluster.  Further details on this work can be found 
in \cite{downs}.

\begin{table}
\begin{center}
\begin{tabular}{|ccccc|}\hline
Sig. dig. &  No. clusters &  Largest cluster &  No. discrep.
 & No. discrep. cl.
 \\ \hline
  1     &      138    &         7037   &           3    &               3 \\
  1     &      148    &         7034   &           1    &               1 \\
  1     &      167    &         6923   &           9    &               7 \\
\hline
\end{tabular}
\end{center}
\caption{Results of k-means for the same three data sets used heretofore,
each relating to 7500 chemical structures, with 1052 descriptors.
``Sig. dig.'': number of significant digits used.  ``No. clusters'': number
of clusters in the data set of 7500 chemical structures, associated with the
number of significant digits used in the Baire scheme.  ``Largest cluster'':
cardinality.  ``No. discrep.'': number of discrepancies found in k-means
clustering outcome.  ``No. discrep. cl.'': number of clusters containing these
discrepant assignments.}
\label{tab4}
\end{table}

\section{Conclusions}

We have been clear in this work in regard to where and when we 
used a Euclidean metric, or other dissimilarity, as input.  
We used correspondence analysis, for
instance, for its property of ``Euclideanizing'' data in the form of counts 
or numbers of occurrences.  Such treatment of the input data was to allow
comparability of results, in a common framework, and in addition it was
noted that very limited assumptions were made in regard to the input data.

It has been our aim in this work to link observed data with an ultrametric
topology for such data.  The traditional approach in data analysis, of course, 
is to impose structure on the data.  This is done, for example, by using 
some agglomerative hierarchical clustering algorithm.  We can always 
do this (modulo distance or other ties in the data).  Then we can assess
the degree of fit of such a (tree or other) structure to our data.  

For our purposes, here, this is unsatisfactory.  

Firstly, our aim was to show 
that ultrametricity can be naturally present in our data, globally or 
locally.  We did not want any ``measuring tool'' such as an 
agglomerative hierarchical clustering algorithm to overly influence 
this finding.  (Unfortunately \cite{ref09b} suffers from precisely this 
unhelpful influence of the ``measuring tool'' of the subdominant 
ultrametric.  In other respects, \cite{ref09b} is a seminal paper.) 

Secondly, let us assume that we did use hierarchical clustering, and then 
based our discussion around the goodness of fit.  This again is a traditional
approach used in data analysis, and in statistical data modeling.  But such 
a discussion would have been unnecessary and futile.  For, after all, if 
we have ultrametric properties in our data then many of the widely used
hierarchical clustering algorithms will give precisely the same outcome, 
and furthermore the fit is by definition exact.  

In linking data with an ultrametric embedding, whether local only, or global,
we have, in this article, 
proceeded also in the direction of exploiting this achievement.  
While some applications, like discrimination between time series signals, 
or texts, have been covered here, other applications like bioinformatics 
database search and discovery, and analysis of large scale 
cosmological structures \cite{murpad}, 
have just been opened up.  In \cite{ezhov} this methodology is applied
to quantum statistics. 
%In the distance there 
%looms the challenge of analysis of networks of enormous size (internet, or 
%biological).  
There is a great deal of work to be accomplished. 

% BibTeX users please use
% \bibliographystyle{}
% \bibliography{}

\begin{thebibliography}{99.}

%\bibitem{alt}
%M.V. Altaisky and B.G. Sidharth, ``p-Adic physics below and above 
%Planck scales'', {\em Chaos, Solitons and Fractals}, {\bf 10}, 167--176, 1999.

\bibitem{bel}
R. Bellman, Adaptive Control Processes: A Guided Tour, 
Princeton University Press, 1961.

\bibitem{ref2}
J.P. Benz\'ecri, L'Analyse des Donn\'ees, Tome I Taxinomie, 
Tome II Correspondances, 2nd ed., Dunod, Paris, 1979.

%\bibitem{bus}
%D. Bustos, G. Navarro  and E. Ch\'avez, 
%Pivot selection techniques for
%proximity searching in metric spaces, Pattern Recognition Letters,
%24, 2357--2366, 2003.

\bibitem{cai1}
F. Cailliez and J.P. Pag\`es, Introduction \`a l'Analyse de
Donn\'ees, SMASH (Soci\'et\'e de Math\'ematiques Appliqu\'ees et de
Sciences Humaines), Paris, 1976.

\bibitem{cai2}
F. Cailliez, The analytical solution of the additive constant
problem, Psychometrika, 48, 305--308, 1983.

%\bibitem{cha1}
%E. Ch\'avez and G. Navarro, Measuring the dimensionality of general
%metric spaces, Technical Report TR/DCC-00-1, Department of Computer
%Science, University of Chile, 2000.

\bibitem{cha2}
E. Ch\'avez, G.  Navarro, R. Baeza-Yates and J.L. Marroqu\'{\i}n, 
Proximity searching in metric spaces, ACM Computing Surveys, 33,
273--321, 2001.

%\bibitem{cha3}
%E. Ch\'avez and G. Navarro, 
%Probabilistic proximity search: fighting
%the curse of dimensionality in metric spaces, Information
%Processing Letters, 85, 39--56, 2003.

\bibitem{desoete}
G. de Soete, 
A least squares algorithm for fitting an ultrametric tree 
to a dissimilarity matrix,
Pattern Recognition Letters, 2, 133--137, 1986.

\bibitem{ezhov}
A.A. Ezhov and A.Yu. Khrennikov, On ultrametricity and a symmetry 
between Bose-Einstein and Fermi-Dirac systems, 
in A.Yu. Khrennikov, Z. Raki\'c and I.V.
Volovich, Eds., p-Adic Mathematical Physics, American Institute of Physics
Conf.\ Proc.\ Vol.\ 826, 55--64, 2006.


\bibitem{fish}
R.A. Fisher, The use of multiple measurements in taxonomic
problems, The Annals of Eugenics,  7, 179--188, 1936.


%\bibitem{frbu1}
%D. Fraix-Burnet, Astrocladistics web page, www-laog.obs.ujf-grenoble.fr/
%$\sim$fraix/astroclad.htm

%\bibitem{frbu2}
% D. Fraix-Burnet, P. Choler, E. Douzery and A. Verhamme
%``Astrocladistics: a phylogenetic analysis of galaxy evolution. 
%I. Character evolutions and galaxy histories'',
%{\em Journal of Classification}, submitted, 2005.

%\bibitem{frbu3}
%D. Fraix-Burnet, E. Douzery, P. Choler and A. Verhamme
%``Astrocladistics: a phylogenetic analysis of galaxy evolution. 
%II. Formation and diversification of galaxies'',  
%{\em Journal of Classification}, submitted, 2005.

%\bibitem{fra}
%C. Fraley and A.E. Raftery, 
%``How many clusters? Which clustering method? Answers via
%model-based cluster analysis'', {\em The Computer Journal}, 
%{\bf 41}, 578--588, 1998.


\bibitem{hall}
P. Hall, J.S. Marron and A. Neeman, 
``Geometric representation of high dimension low sample size data'',
{\em Journal of the Royal Statistical Society B}, {\bf 67}, 427--444, 
2005.

\bibitem{hornik}
K. Hornik, 
A CLUE for CLUster Ensembles, 
Journal of Statistical Software, 14 (12), 2005.

\bibitem{ref07}
I.C. Lerman, Classification et Analyse Ordinale des Donn\'ees, 
Paris, Dunod, 1981.

%\bibitem{moulder}
%J. Moulder,  Electromagnetic fields and human health: Power lines
%and cancer FAQs, http://www.mcw.edu/gcrc/cop/powerlines-cancer-FAQ/toc.html,
%2006.

\bibitem{mur85a}
F. Murtagh, Multidimensional Clustering Algorithms, Physica-Verlag, 1985.

%\bibitem{mur85}
%F. Murtagh and A. Heck, {\em Multivariate Data Analysis}, Kluwer, 1987.

\bibitem{ref08} 
F. Murtagh, 
On ultrametricity, data coding, and computation, 
Journal of Classification, 21, 167--184, 2004.

\bibitem{ref088}
F. Murtagh, Identifying the ultrametricity of time series, 
European Physical Journal B, 43, 573--579, 2005.

\bibitem{ref0888}
F. Murtagh, A note on local ultrametricity in text,
http://arxiv.org/pdf/cs.CL/0701181, 2007.
 
\bibitem{ref08888}
F. Murtagh, Correspondence Analysis and Data Coding with R and Java,
Chapman \& Hall/CRC, 2005.

\bibitem{murpad}
F. Murtagh, From data to the physics using ultrametrics: new results
in high dimensional data analysis, 
in A.Yu. Khrennikov, Z. Raki\'c and I.V.
Volovich, Eds., p-Adic Mathematical Physics, American Institute of Physics
Conf.\ Proc.\ Vol.\ 826,
151--161, 2006.

\bibitem{downs}
F. Murtagh, G. Downs and P. Contreras, ``Hierarchical clustering of 
massive, high
dimensional data sets by exploiting ultrametric embedding'', 2007, submitted.  

\bibitem{autoonto}
F Murtagh, J. Mothe and K. Englmeier, ``Ontology from local hierarchical
structure in text'', http://arxiv.org/abs/cs.IR/0701180, 2007.

%\bibitem{newman}
%M.E.J. Newman, Power laws, Pareto distributions and Zipf's law,
%Contemporary Physics, 46, 323--351, 2005.
%%cond-mat/0412004
%%http://arxiv.org/abs/cond-mat/0412004

\bibitem{neu}
E. Neuwirth and L. Reisinger,  
Dissimilarity and distance coefficients
in automation-supported thesauri,  Information Systems, 
7, 47--52, 1982.

%\bibitem{oneill}
%E. O'Neill, Understanding ubiquitous computing: a view from
%HCI, in Discussion following R. Milner, Ubiquitous computing: 
%how will we understand it?, Computer Journal, in press, 2006.

\bibitem{ref09a}
R. Rammal, J.C. Angles d'Auriac  and B. Doucot, 
On the degree of ultrametricity, 
Le Journal de Physique -- Lettres, 46, L-945--L-952, 1985.

\bibitem{ref09b}
R. Rammal, G. Toulouse and M.A. Virasoro, 
Ultrametricity for physicists, 
Reviews of Modern Physics, 58, 765--788, 1986.
% 765--788

%\bibitem{schmid}
%H. Schmid,  Probabilistic part-of-speech tagging using decision 
%trees, Proc. Intl. Conf. New Methods in Language Processing.  
%TreeTagger, www.ims.uni-stuttgart.de/projekte/corplex/TreeTagger/
%DecisionTreeTagger.html, 1994.

\bibitem{tor} 
W.S. Torgerson, Theory and Methods of Scaling, Wiley, 1958.

\bibitem{ref120}
A. Treves, 
On the perceptual structure of face space,
BioSystems, 40, 189--196, 1997.
%189--196.

%\bibitem{vanr}
%C.J. van Rijsbergen, Information Retrieval, 2nd ed.
%Butterworths, 1979.

%\bibitem{wstalk}
%WS-Talk: Web Services Communicating 
%in the Language of Their User Community, Sixth Framework CRAFT 
%(Small and Medium Enterprises Cooperative Research) project, 2005--2006, \\
%http://thames.cs.rhul.ac.uk/$\sim$wstalk.

\end{thebibliography}
%
% Non-BibTeX users please follow the syntax
% the syntax of ``referenc.tex'' for your own citations
%\input{referenc}

\end{document}